\input amstex
\documentstyle{amsppt}
\input amssym
\input psfig.sty 
\magnification=\magstep1
\def\tvi{\vrule height 12truept depth 5truept width0truept}
\def\tv{\tvi\vrule}
\def\centro#1{\kern .7em\hfill#1\hfill\kern .7em}
\def\hline{\noalign{\hrule}}

\newcount\escala 
\newbox\texto
\newcount\dibno
\dibno=1 
\def\ilustracion#1,#2(#3)#4.{
\centerline{\psfig{figure=#3.ps,width=#2mm,height=#1mm}}%
\nobreak{\setbox\texto=\hbox%
{\it Fig. \the\dibno. #4}
\ifdim\wd\texto>\hsize{\narrower\noindent\unhcopy\texto\par}%
\else\centerline{\copy\texto}\fi}
\global\advance\dibno by 1
\ifdim\lastskip<\smallskipamount
\removelastskip\penalty55\medskip\fi
  \ignorespaces\smallbreak}
  
\newcount\parno \parno=1
\newcount\sparno \sparno=1
\newcount\prono \prono=1
\def\sec{\S\the\parno.-\ \global\sparno=1 \global\prono=1}
\def\ssec{\the\parno.\the\sparno-\ }
\def\etiqueta{\hbox{(\the\parno.\the\prono)}}
\def\finparrafo{\global\advance\parno by1
\vskip.1truecm\ignorespaces}
\def\finsparrafo{\global\advance\sparno by1
\vskip.1truecm\ignorespaces}
\def\finparrafo{\global\advance\parno by1
\vskip.1truecm\ignorespaces}
\def\cita{\ignorespaces\ \the\parno.\the\prono%
\global\advance\prono by 1}

\def\bc{{\Bbb C}} \def\ba{{\Bbb A}}
\def\bp{{\Bbb P}} \def\bn{{\Bbb N}}
 \def\bz{{\Bbb Z}}
\def\br{{\Bbb R}}

\def\F{{\Cal F}}
\def\U{{\Cal U}}
\def\T{{\Cal T}}
\def\C{{\Cal C}}
\def\Caf{{\Cal C}_{\text{af}}}
\def\Ccaf{{\Cal C}_{1,\text{af}}}
\def\Cqaf{{\Cal C}_{2,\text{af}}}
\def\L{\ell}
\def\Li{\L_{\infty}}

\def\FF{{\Bbb F}}

\def\ba{{\Bbb A}}
\def\bz{{\Bbb Z}}
\def\bn{{\Bbb N}}
\def\AA{{\Bbb A}}
\def\PP{{\Bbb P}}
\def\QQ{{\Bbb Q}}

\def\CC{{\Bbb C}}

\def\ZZ{{\Bbb Z}}

\def\cL{\mathop {\Cal L}\nolimits}

\def\Sing{\text{Sing}}
\def\mod{\text{mod}}

\def\NS{\text{NS}}

\def\Pic{\text{Pic}}
\def\Supp{\text{Supp}}
\def\hiya#1{\mathrel{\mathop{\leftarrow}\limits^#1}}
\def\miya#1{\mathrel{\mathop{\rightarrow}\limits^#1}}

\def\D{{\Cal D}}
\def\disc{\text{disc}}
\def\tpi {\tilde{\pi}}

\def\tS{\tilde S}

\def\qed{\quad $\square$}

\NoBlackBoxes

\refstyle{A}
\widestnumber\key{A et al.} 
\NoRunningHeads  

\topmatter

\title
Sextics with singular points in special position
\endtitle

\abstract
In this paper we show a Zariski pair of sextics which is 
not a degeneration of the original example given by Zariski.
This is the first example of this kind known. 
The two curves of the pair have a trivial
Alexander polynomial. The difference in the topology of 
their complements can only be detected via finer invariants
or techniques. In our case the generic braid monodromies, 
the fundamental groups, the characteristic varieties and the 
existence of dihedral coverings of $\PP^2$ ramified along them
can be used to distinguish the two sextics. Our intention is not 
only to use different methods and give a general description 
of them but also to bring together different perspectives of 
the same problem.
\endabstract

\author
E. ARTAL BARTOLO$^1$, J. CARMONA RUBER$^1$, J.I~COGOLLUDO$^2$
and Hiro-o TOKUNAGA$^3$
\endauthor

\leftheadtext{E. ARTAL, J. CARMONA, J.I~COGOLLUDO
AND H. TOKUNAGA}

\address
Departamento de Matem\'aticas,
Universidad de Zaragoza,
Campus Plaza San Francisco s/n
E-50009 Zaragoza SPAIN
\endaddress
\email
artal\@posta.unizar.es
\endemail

\address
Departamento de Sistemas inform\'aticos y programaci\'on,
Universidad Compluten\-se,
Ciudad Universitaria s/n
E-28040 Madrid SPAIN
\endaddress
\email
jcarmona\@eucmos.sim.ucm.es
\endemail

\address
Departamento de \'Algebra,
Universidad Complutense,
Ciudad Universitaria s/n
E-28040 Madrid SPAIN
\endaddress
\email
jicogo\@eucmos.sim.ucm.es
\endemail

\address
Department of Mathematics,  
Tokyo Metropolitan University,
1-1 Minamiohsawa
Hachioji 192-0397
Tokyo JAPAN
\endaddress
\email
tokunaga\@comp.metro-u.ac.jp
\endemail

\thanks
$^1$ Partially supported by
DGES PB97-0284-C02-02
\endthanks

\thanks
$^2$ Partially supported by
DGES PB97-0284-C02-01
\endthanks

\thanks
$^3$ Partially supported by 
JSPS 11640034
\endthanks


\endtopmatter

\document

The starting point of this paper is the existence of 
a certain pair of sextic curves 
$\C^{(1)}, \C^{(2)} \subset \PP^2$ in the complex 
projective plane. These curves have the same combinatorial 
properties and hence are candidates to form a Zariski pair. 
Roughly speaking, two curves form a Zariski pair if they have 
the same combinatorial data (degree of each irreducible 
component, local type of singularities,...) but different 
embeddings in $\PP^2$ --\,cf. \cite{A}. The combinatorial 
data of our pair are the following: both are reducible 
sextics, given by a smooth conic and a quartic; the quartic
has two singular points, say $P_1$ and $P_2$, of types $\AA_1$ 
and $\AA_3$ respectively; the conic and the quartic intersect 
at one single point $Q$ of type $\AA_{15}$. The non-combinatorial 
datum that distinguishes them is that the tangent line at $Q$
contains (for $\C^{(2)}$) or not (for $\C^{(1)}$) the point $P_2$. 
The easiest and most common strategy to distinguish 
embeddings is to compute the Alexander polynomial of the curves. 
Such a polynomial is known to be sensitive to the special 
position of some singular points of the curves 
--\,cf. \cite{A}, \cite{L1}, \cite{D}. 
In our case, even though the position of singularities 
distinguishes the two embeddings, it has no effect on the 
Alexander polynomial, which coincides in both cases. 
Examples of Alexander-equivalent Zariski pairs are already 
known --\,cf. \cite{A-Ca}, \cite{O}. 
In this paper we will describe four different techniques to 
prove that two given curves constitute a Zariski pair and will
apply them to the sextics $\C^{(1)}$ and $\C^{(2)}$:

\smallbreak\item{--} Computation of fundamental groups of the
complement. Let $\C$ be an algebraic curve in the complex 
projective plane $\PP^2$. The fundamental group of the complement
$\PP^2 \setminus \C$ is a topological invariant of the pair
$(\PP^2,\C)$. Computation of such groups is in general a 
difficult task, and even if one manages to compute them, by 
finding a finite presentation, it might not be possible to decide 
whether or not the groups are isomorphic. In our case we use a 
non-generic braid monodromy to recover a presentation of the 
fundamental group.

\smallbreak\item{--} Computation of generic braid monodromy 
groups. This is also a topological invariant of the pair
$(\PP^2,\C)$, but it is finer than the fundamental group. 
In fact, it determines the homotopy type of the complement 
$\PP^2 \setminus \C$ --\,cf. \cite{L2}. 
It might be useful in cases where the fundamental groups are 
either isomorphic or at least difficult to compare. We think that
it is worth looking for effective invariants, e.g. in the aim
of the work of Libgober --\,cf. \cite{L3}.

\smallbreak\item{--} Computation of characteristic varieties.
First defined by Libgober \cite{L4}, they are a generalization 
of Alexander polynomials. If Alexander polynomials are related 
to (infinite) cyclic coverings of $\PP^2$ ramified along a curve 
$\C$, the characteristic varieties are related to general abelian 
coverings of such kind. They might be computed from the fundamental 
groups; in this case they provide effective invariants. It is also 
possible to compute them without knowing the fundamental group 
using the theory of ideals of quasiadjunction, which generalizes 
computations of Alexander polynomials. This technique can be followed 
in detail in \cite{L4}.

\smallbreak\item{--} Existence of dihedral coverings. The work of 
the fourth author has proven that it is also an effective way to 
obtain properties of the complements without an explicit 
computation of fundamental groups. The existence of a certain
kind of dihedral coverings of $\PP^2$ branched along $\C$ can be 
characterize by the existence of divisors satisfying certain 
algebraic properties. This allows to establish a connection between 
the position of singularities and the existence of non-abelian
coverings of $\PP^2$ ramified along $\C$.

\bigbreak
As mentioned in the abstract, this example is the first known
Zariski pair of sextics that is not a degeneration of the 
original example given by Zariski in \cite{Z}. Namely, there 
is no conic passing through six --\,maybe infinitely near\,-- 
non-nodal singular points of either sextic. 

\smallbreak
The authors have found other pairs that might be worth 
studying. We won't discuss them here due to length considerations,
but we include a brief description of them. The combinatorial 
data are the following: reducible sextics whose components 
are a quintic and a line; the quintic has two singular points, 
say $P_1$ and $P_2$, of types $\AA_1$ and $\AA_9$ respectively; 
the quintic and the line intersect at a single point $Q$ of type 
$\AA_9$. There are three such curves up to projective change 
of coordinates. One of them, say $\D^{(2)}$, has the property 
that the tangent line at $P_2$ contains the point $Q$. This curve 
admits a birrational transformation to the sextic $\C^{(2)}$ 
considered in this paper. The other two sextics, $\D^{(1,1)}$ 
and $\D^{(1,2)}$, are conjugated in $\QQ(\sqrt{5})$ and don't 
admit a birrational transformation to $\C^{(1)}$. Both pairs
$(\D^{(1,1)},\D^{(2)})$ and $(\D^{(1,2)},\D^{(2)})$ are
Zariski pairs and can be distinguished using the techniques
presented in this paper, whereas the pair 
$(\D^{(1,1)},\D^{(1,2)})$ remains undistinguishable for us.

\bigbreak
A summarized description of the sections could be the following.
In \S 1 we define and construct the curves. In section \S 2 we
briefly describe the braid monodromy of affine curves
and the Zariski-van Kampen method as well as the connection 
between fundamental groups of affine and projective curves.
Sections \S 3 and \S 4 are devoted to compute non-generic
braid monodromies for $\C^{(1)}$ and $\C^{(2)}$, derive
presentations for the fundamental groups of their complements
and show that they are non-isomorphic. In Section \S 5 we 
compute the generic braid monodromies of $\C^{(1)}$ and 
$\C^{(2)}$. In section \S 6 we introduce the characteristic
varieties and briefly describe the procedure to calculate
them in general following \cite{L4}. We apply such a method
to our pair of sextics showing that both differ in the
characteristic variety of depth one. The last section, \S 7, 
is devoted to describing the connection between a special
type of dihedral coverings of $\PP^2$ ramified along a curve 
and the existence of a divisor satisfying certain algebraic 
properties in terms of the Picard group. We prove that there 
exists such a dihedral covering of degree 16 ramified along 
one sextic whereas there is none for the other.

\head\sec Definition of the curves
\endhead

Let $\C$ be an algebraic curve in the complex projective 
plane $\PP^2$ satisfying the following properties:

\smallbreak\item{P1.} $\C$ has two irreducible components 
$\C_1$ and $\C_2$.

\smallbreak\item{P2.} $\C_1$ is a smooth conic.

\smallbreak\item{P3.} $\C_2$ is a curve of degree four having two 
singular points of type $\ba_3$ (tacnode) and $\ba_1$ (node).

\smallbreak\item{P4.} $\C_1\cap\C_2=\{P\}$, and 
$P$ is a smooth point of $\C_2$.
Then $(\C,P)$ is a singularity of type $\ba_{15}$.

\newbox\ejemi
\global\setbox\ejemi=\etiqueta
\definition{Example\cita} Let us consider the curve 
$\C^{(1)}=\C^{(1)}_1 \cup \C^{(1)}_2$ where:
$$ \C^{(1)}_1: 8 z^2 -20 y z+36 x z+17 y^2-18 x y=0,$$
$$ \C^{(1)}_2: 8 x^2 z^2 -16 x y^2 z+52 x^2 y z -36 x^3 z- 
37 x^2 y^2+18 x^3 y-y^4+20 x y^3=0.$$
This curve verifies the conditions P1-P4, where $[1:0:0]$ 
is the singular point $P$ of type $\ba_{15}$, $[0:0:1]$ is 
of type $\ba_3$ and $[1:1:0]$ is of type $\ba_1$. Note 
that the common tangent line at $P$ is $y=2z$ and 
it does not pass through $\ba_3$.
\enddefinition

\vskip 10pt

\ilustracion 36,46(curvag)Affine curve in $(y,z)$.

\newbox\ejemii
\global\setbox\ejemii=\etiqueta
\definition{Example\cita} Let us consider the curve 
$\C^{(2)}=\C^{(2)}_1 \cup \C^{(2)}_2$ where:
$$\C^{(2)}_1: 3 x^2+2 x y+ 108 z^2=0,$$
$$\C^{(2)}_2:2 x y^3+3 x^2 y^2+ 108 y^2 z^2 -x^4=0.$$
This curve also verifies P1-P4 above, where $[0:1:0]$ is 
the singular point $P$ of type $\ba_{15}$, $[0:0:1]$ is 
of type $\ba_3$ and $[1:-1:0]$ is of type $\ba_1$. Note 
that the common tangent line at $P$ is $x=0$ and it 
passes through $\ba_3$.
\enddefinition

\vskip 10pt

\ilustracion 36,46(curvang)Affine curve in $(x,z)$.

The way these examples can be constructed will be useful
in the future. Let us suppose that $\C$ satisfies P1-P4.
We want to find a suitable Cremona transformation to $\PP^2$
that simplifies the problem. Let us consider the two dimensional
family $\bp$ of conics passing through $\ba_3$ and $\ba_{15}$ 
which are tangent to $\C$ at $\ba_3$. There is a birational map 
$\bp^2\dasharrow\check\bp$, where $\check\bp$ is the dual of 
$\bp$. Let us identify $\bp^2$ with $\check\bp$. There is an 
easy description of the mapping in terms of blowing-ups and 
blowing-downs as follows:

\smallbreak\item{--} Let us denote by $T$ the tangent line to $\C$ 
at $\ba_3$ and by $L$ the line through $\ba_{15}$ and $\ba_3$.
Blow up $\ba_{15}$ and $\ba_3$. Let us denote $A_{15}$ and $A_3$
the corresponding exceptional components. Let us denote by $A$ the
intersection of $T$ and $A_3$, and by $P'$ the intersection of the 
proper transform of $\C$ and $A_{15}$.

\smallbreak\item{--} Let us blow up $A$ and blow down $L$. The 
corresponding exceptional component will be also denoted by $A$.

\smallbreak\item{--} Let us blow down $A_3$ and $T$ and denote 
accordingly the contracted points.
\smallbreak
The above family of conics has been transformed into the family of 
conics passing through $T$ and tangent to $A_{15}$ at $A_3$.

It is an easy exercise to show that $\C_2$ has been transformed onto 
a nodal cubic (nothing happened near $\ba_1$) transversal to $A$ and
tangent to $A_{15}$ at $A_3$. Moreover, $A_3$ is an inflection point
of this cubic if and only if $\C$ is tangent to $L$ at $\ba_{15}$.
On the other hand, the transform of $\C_1$ is an irreducible cubic 
also tangent to $A_{15}$ at $A_3$ and having a double point at $T$.
There are two possible cases for the intersection of the transforms
of $\C_1$ and $\C_2$:

\smallbreak\item{(a)} If $\C$ is not tangent to $L$ at $\ba_{15}$,
then the transforms are ordinary tangents at $A_3$ and the contact 
order at $P'$ is equal to seven.

\smallbreak\item{(b)} If $\C$ is tangent to $L$ at $\ba_{15}$,
then the contact order at $A_3$ is equal to nine.

\smallbreak
Heuristically the second case is a degeneration of the first one 
and they may be treated together. We will also denote by $\C_1$ 
and $\C_2$ the transforms of the given curves by the above Cremona 
transformation.

\newbox\ellipt
\global\setbox\ellipt=\etiqueta

\proclaim{Proposition\cita} Up to projective transformation
there are only two pairs of curves matching either (a) or (b)
above. Namely:
\smallbreak\item{\text(1)}
There exists a pair of nodal cubics $\C^{(1)}_1$ and $\C^{(1)}_2$ 
satisfying:
\smallbreak\itemitem{\text(a1)}
$\C^{(1)}_1 \cap \C^{(1)}_2=\{P_1,P'\}$ where 
$P_1$ and $P'$ are singular points of 
$\C^{(1)}=\C^{(1)}_1\cup\C^{(1)}_2$ of type $\ba_{11}$ and $\ba_3$ 
respectively, and
\smallbreak\itemitem{(b1)}
the tangent line to $\C^{(1)}$ at $P'$ contains the point $P_1$.

These cubics generate a pencil containing a reducible curve having 
as components the common tangent line at $P'$ and the conic which 
has intersection number six at $P_1$ with both cubics.

\smallbreak\item{\text (2)}
There exists a pair of nodal cubics $\C^{(2)}_1$ and $\C^{(2)}_2$ 
such that $\C^{(2)}_1 \cap \C^{(2)}_2=\{P\}$ is a singular point of 
$\C^{(2)}=\C^{(2)}_1\cup\C^{(2)}_2$ of type $\ba_{17}$.

These cubics generate a pencil containing a triple line which is
the common tangent line at $P$, an inflection point for both curves.
\endproclaim

\demo{Proof} Let $\C_2$ be any nodal cubic. We will denote by
$\varphi: \bc^* \rightarrow Reg(\C_2)$ any parametrization
inducing an isomorphism of groups. The group structure of 
$Reg(\C_2)$ is the geometrical one where the zero element is 
an inflection point. Such a parametrization is well defined up 
to a different choice of inflection point 
(obtained by multiplying by a cubic root of unity) and 
interchanging the branches of the nodal point (obtained by 
taking the inverse). In other words, there is a transitive and 
free action by the dihedral group $\D_6$ of order $6$ on the 
space of all such parametrizations.

We recall that the geometrical group structure can be described 
as follows: let $P_1,\dots,P_r\in Reg(\C_2)$ and let 
$m_1,\dots,m_r\in\bn$ such that $m_1+\dots+m_r=3n$. Then 
$m_1P_1+\dots+m_rP_r=0$ if and only if there exists a curve $D$ 
of degree $n$ such that $(D\cdot\C_2)_{P_j}=m_j$, $j=1,\dots,k$.

Let us denote $t_1=\varphi^{-1}(P')$ and $t_2=\varphi^{-1}(A_3)$. 
Then:
$$t_1 t_2^2=1,\qquad t_1^7 t_2^2 =1,$$
i.e., $t_1=t_2^{-2}$ and $t_2^{12}=1$. The factorization 
$t^{12}-1=(t^3-1)(t^3+1)(t^6+1)$ represents the space of 
orbits of the roots of $t^{12}-1$ by the action of $\D_6$.

Let $t_2$ be a $12$th root of unity. Consider the curve $\C_2$ 
and consider the pencil of cubics having contact order equal 
to seven (resp. two) with $\C_2$  at $P'$ (resp. $A_3$).
To prove the existence of $\C_1$ we need to find a nodal cubic
in this pencil (whose node is not $A_3$).

One has several possibilities for $t_2$:

\smallbreak\item{--} If $t_2^3=-1$, then $t_1^3=1$. Then
$P'$ is an inflection point of $\C_2$. As it is readily
seen, the double tangent to $P'$ and the line $A_{15}$ 
joining $P'$ and $A_3$ are a member of the pencil. 
Blowing up the base points results in the elliptic fibration 
associated to the pencil. An easy argument of euler 
characteristics shows that there is no such a fiber as $\C_1$.

\smallbreak\item{--} If $t_2^6=-1$, one obtains (1) as follows.
Denote
$P'=\varphi(t_2)$. Note that
$$t_1^6=(-t_2^4)^6=1.$$
Moreover, $t_1$ has order six. Hence, the tangent line $L'$ to 
$\C^{(1)}_2$ at $P'$ contains the point $\varphi(t_1)=P_1$. 
In fact, there is smooth conic $Q$ intersecting $\C^{(1)}_1$ 
only at $P_1$ (i.e. with multiplicity six). The cubics
$\C'=Q + L'$ and $\C^{(1)}_2$ define a pencil and, after
blowing up the base points, an elliptic fibration from a
rational surface $X_{8211}$ on $\PP^1$. This fibration has
three evident exceptional fibers: the one given by $\C^{(1)}_2$ 
($I_1$ in Kodaira's notation), the one given by $\C'$ ($I_8$)
and a nodal cubic with the node on $P'$ ($I_2$). Since
the euler characteristic of $X_{8211}$ is $12$ and
$\chi(I_n)=n$, there must be another special fiber $I_1$.
This fiber corresponds to a nodal cubic $\C^{(1)}_1$, whose 
node is not a base point of the pencil. Hence, this is the
required pair of cubics.

\smallbreak\item{--} If $t_2^3=1$, one obtains (2) as follows.
In this case $P'=\varphi(t_2)$ is an inflexion point of 
$\C^{(2)}_2$. Let $L'$ be the tangent line of $\C^{(2)}_2$
at $P$. Consider the pencil defined by $\C^{(2)}_2$ and
$3L'$. As in the previous case, after blowing up the nine 
base points of the pencil (counted with multiplicity), we 
obtain a rational surface $X_{211}$ defining an elliptic 
fibration over $\PP^1$. In this case, we have two evident special 
fibers, namely the one determined by $\C^{(2)}_2$ ($I_1$) 
and the one determined by the triple line ($II^*$). Since 
$\chi(II^*)=10$, $\chi(I_1)=1$ and $\chi(X_{211})=12$, 
there has to be another special fiber whose euler 
characteristic equals one, that is, a nodal cubic $\C^{(2)}_1$ 
whose node is not on $P$. The pair $\C^{(2)}_1$, $\C^{(2)}_2$ 
is the required pair of cubics.\qed
\enddemo

With a suitable choice of equations for $\C_i^{(j)}$ 
and for the Cremona transformations one can obtain the 
equations given in \copy\ejemi\ and \copy\ejemii.

\proclaim{Remark \cita} The existence of the elliptic surfaces
$X_{8211}$ and $X_{211}$ can also be derived directly from 
\cite{Mi-P1}. Both rational surfaces will be specially useful 
in section \S 7.
\endproclaim

Along these sections we will prove the following

\newbox\zarpair
\global\setbox\zarpair=\etiqueta
\proclaim{Main theorem\cita}
The sextic curves $\C^{(1)}$ and $\C^{(2)}$ form
a Zariski pair.
\endproclaim

\finparrafo

\head\sec Some facts about fundamental groups and 
braid monodromy of affine curves
\endhead

In the first part of this section we will define the
basic concepts required to work with braid monodromies
of affine curves. In the second part we will prove the
relationship between fundamental groups of projective
and affine curves.

\medbreak
Let $a_1,\dots,a_k$ be $k$ distinct complex numbers and
put $\U=\bc\setminus\{a_1,\dots,a_k\}$. We may suppose that, 
for $1\leq i<k$ we have either $\Re a_i> \Re a_{i+1}$ or
$\Re a_i= \Re a_{i+1}$ and $\Im a_i> \Im a_{i+1}$. Let us 
choose $a_0\in\br$ such that $a_0\gg\max\{|a_1|,\dots,|a_k|\}$. 
One can construct an oriented piecewise linear curve $\Gamma$
which is the union of the segments $[a_{j-1},a_j]$. Take 
$0<\varepsilon\ll 1$ such that the closed disks $D_i$ 
of radius $\varepsilon$ centered at $a_i$ , $i=1,\dots,k,$ 
are pairwise disjoint and do not contain $a_0$. 
Let us orient $\alpha_i:=\partial D_i$ counterclockwise.
Let us denote by $p_i^+$ the first point of intersection 
of $\Gamma$ and $\alpha_i$. The remaining point of intersection
of $\Gamma$ and $\alpha_i$ will be denoted by $p_i^-$.
The path $\Gamma$ cuts $\alpha$ into two components $\alpha_i^+$
and $\alpha_i^-$, where $\alpha_i^+$ starts at $p_i^+$. 
We also construct $\Gamma_j$, $j=1,\dots,k$ as follows:

\smallbreak\item{--} $\Gamma_1$ is the subpath in $\Gamma$ from
$a_0$ to $p_1^+$.

\smallbreak\item{--} If $j>1$ then $\Gamma_j$ is the subpath in 
$\Gamma$ from $p_{j-1}^-$ to $p_j^+$.

\smallbreak
For $j=1,\dots,s$ let 
$\beta_j:=\Gamma_1\cdot\prod_{i=2}^j(\alpha_i^+\cdot\Gamma_i)$,
which is a path from $a_0$ to $p_j^+$. We define:
$$\gamma_j:=\beta_j\cdot\alpha_j\cdot\beta_j^{-1}.$$
It is well known that the homotopy classes of
$\gamma_1,\dots,\gamma_k$ form a basis of the free group
$\pi_1(\U,a_0)$.

\definition{Definition\cita} An ordered family of
lassos $\mu_1,\dots,\mu_k$ is called
a {\bf standard basis} if these paths are
homotopy equivalent to $\gamma_1,\dots,\gamma_k$.
\enddefinition

\definition{Definition\cita} Let $a_1,\dots,a_k$ be as 
above and let $b_1,\dots,b_k$ be another set of $k$ distinct 
points in $\bc$. A {\bf motion} from $a_1,\dots,a_k$ to
$b_1,\dots,b_k$ is a set of paths
$\delta_i\:[0,1]\to\bc$, $i=1,\dots,k$ such that:

\smallbreak\item{--} $\delta_i(0)=a_i$, $\delta_i(1)=b_{\rho(i)}$,
$i=1,\dots,s$, where $\rho$ is a permutation of $\{1,\dots,s\}$.

\smallbreak\item{--} For all $t\in[0,1]$, the set 
$\{\delta_i(t)\}_{i=1}^s$ contains $s$ distinct points.

If $a_i=b_i$, $i=1,\dots,k$, we say that the motion is based
on $a_1,\dots,a_k$.
\enddefinition

\remark{Remark} Note that the order of the $k$ points is 
not relevant.
\endremark

\definition{Definition\cita} Two motions
$\delta_i^0\:[0,1]\to\bc$ and
$\delta_i^1\:[0,1]\to\bc$, $i=1,\dots,s$ from on 
$a_1,\dots,a_s$ to $b_1,\dots,b_k$  
are homotopic if there exists a set of homotopies 
$H_i\:[0,1]\times [0,1]\to\bc$, $i=1,\dots,k$,
such that: 
\smallbreak\item{--} $\delta_i^0=H_i(\text{--},0)$, 
$\delta_i^1=H_i(\text{--},1)$, $\{1,\dots,k\}$, and

\smallbreak\item{--} if we define $\delta_i^s:=H_i(\text{--},s)$,
then $\{\delta_i^s\}_{i=1}^k$ is also a motion from $a_1,\dots,a_k$
to $b_1,\dots,b_k$.

A braid from $a_1,\dots,a_k$
to $b_1,\dots,b_k$ (resp. a braid based on $a_1,\dots,a_k$) 
is an equivalence class of motions
by this relation. The set of braids based on $a_1,\dots,a_k$
is denoted
$B_{a_1,\dots,a_k}$ and is a group which is of course isomorphic
to the braid group $B_k$.
\enddefinition

The canonical presentation of $B_k$ is given by
$$
\left|\sigma_1,\dots,\sigma_{k-1} : [\sigma_i,\sigma_j]=1\text{ if }
|i-j|\geq 2, \sigma_i\sigma_{i+1}\sigma_i=\sigma_{i+1}\sigma_i\sigma_{i+1},
i=1,\dots,k-2\right|.
$$
We can define a model for these generators in our group. Take 
$1\leq i<k$. Consider the segment $[a_i,a_{i+1}]$ and take a 
small (topological) disk such that this segment is a diameter. Then:

\smallbreak\item{--} Move from $a_i$ to
$a_{i+1}$ along one of the arcs of the boundary of this disk
counterclockwise;

\smallbreak\item{--} move from $a_{i+1}$ to
$a_i$ along the other arc; 

\smallbreak\item{--} take constant paths for the other points.

\vskip 10pt

\ilustracion 30,100(trenza0) Motion between $a_i$ and $a_{i+1}$.

There are several equivalent definitions of the action
of $B_{a_1,\dots,a_k}$ on the free group $\pi_1(\U,a_0)$.
In our case we have for $i=1,\dots,k-1$:
$$\matrix
\gamma_i^{\sigma_i}&=&\gamma_{i+1},\hfill&\\
\gamma_{i+1}^{\sigma_i}&=&\gamma_{i+1}\gamma_i\gamma_{i+1}^{-1},&\\
\gamma_j^{\sigma_i}&=&\gamma_j,&\text{ if }j\neq i,i+1.\\
\endmatrix$$

Now let us recall what is meant by braid monodromy for affine 
curves. Let $\Caf:=\{(x,y)\in\bc^2\mid f(x,y)=0\}$, where
$f(x,y)\in\bc[x,y]$ is of the form:
$$f(x,y)=y^d+a_1(x) y^{d-1}+\dots+a_{d-1}(x)y+a_d(x),
\qquad a_j\in\bc[x],$$
and such that no vertical line is contained in $\Caf$.

Let $S:=\{x_1,\dots,x_r\}$, the set of zeros of $\Delta_y(x)$, 
which is the discriminant of $f$ with respect to $y$. The 
lines $x=x_j$, $j=1,\dots,r$ are the vertical lines which are 
not transversal to $\Caf$, i.e., such that the number
of points of intersection with $\Caf$ is less than $d$. 
Let $x_0\in\bc\setminus S$ such that $x_0\in\br$ and 
$x_0\gg\max\{|x_1|,\dots,|x_r|\}$. Let $y_1,\dots,y_d$ the 
$d$ distinct roots of $f(x_0,y)=0$. We will use the 
following notation:

\smallbreak\item{--} $\F$ for the vertical line $x=x_0$;

\smallbreak\item{--} $\check \F:=\F\setminus \Caf$, i.e., 
$\check \F=\{(x_0,y)\in\bc^2\mid y\neq y_j, j=1,\dots,d\}$.
In order to have a suitable base point we choose $y_0\in\br$
such that $y_0\gg\max\{|y_1|,\dots,|y_d|\}$.

\smallbreak\item{--} $B_d$ for the braid group based on 
$y_1,\dots,y_d$.

Since outside $S$ the curve defines a multi-valuated function,
lassos based on $x_0$ define motions based on $y_1,\dots,y_d$,
which respect homotopy equivalences. This defines a homomorphism
$$\Phi\:\pi_1(\bc\setminus S)\to B_d,$$ 
which is called the {\bf braid monodromy of the curve $\Caf$ 
with equation $f(x,y)=0$} --\,cf. \cite{Mo}. 

Let us recall the Zariski-Van Kampen theorem with this notation. 
Let us consider the group $F:=\pi_1(\check \F,(x_0,y_0))$ which 
is a free group on $d$ generators. Fix a standard basis 
$\mu_1,\dots,\mu_d$ and consider the action of $B_d$. Then,
$$\pi_1(\bc^2\setminus\Caf;(x_0,y_0))=|a\in\F: 
a=a^{\Phi(b)}, a\in\F,b\in \pi_1(\bc\setminus S,x_0)|.$$
Taking a standard basis $\varphi_1,\dots,\varphi_r$ of 
$\pi_1(\bc\setminus S)$ one obtains a finite presentation of this 
group as follows:
$$\pi_1(\bc^2\setminus\Caf;(x_0,y_0))=|\mu_1,\dots,\mu_d: 
\mu_i=\mu_i^{\Phi(\varphi_j)}, i=1,\dots,d,\ j=1,\dots,r|.$$

\remark{Remark\cita} Let us consider a path
$\beta\:[0,1]\to\bc\setminus S$ and let
us denote $R_t$ the set of roots of $f(\beta(t),y)$. 
This path defines a motion from $R_0$ to $R_1$.
The conventions we have followed in the construction
of standard basis allow us to construct well-defined
braids from $y_1,\dots,y_d$ to any subset of $d$ distinct 
points of $\bc$. Then joining $y_1,\dots,y_d$ to $R_0$ and 
$R_1$ one can associate to every $\gamma$ a braid in $B_d$.
If a lasso $\gamma$ is presented as a product of non-closed
paths $\beta_1\cdot\dots\cdot\beta_l$, the braid associated 
to $\gamma$ is the product of the braids associated to each 
$\beta_i$.
\endremark

One can interpret motions of $k$ points as sets of $k$
non-intersecting curves in $\bc\times[0,1]$ such that
each of these curves cut each horizontal plane $\bc\times\{t\}$
in one point, $t\in[0,1]$. In order to draw effectively
this three dimensional picture, one chooses a projection
$\bc\to\br$ such that: $(i)$ at most two points of the 
motion have the same intersection point and $(ii)$ the set 
of such points is isolated; the {\it visible\/} part
of $\bc$ is the negative imaginary semiplane. We will use 
the following criteria in order to choose the projection:

\smallbreak\item{--} Generically the projection $\bc\to\br$
is defined by $z\mapsto\Re z$. If two points have the same 
real part, we mark as visible the one with smaller imaginary 
part.

\smallbreak\item{--} Whenever this projection is not generic
(for example if one has couples of conjugate non-real points),
we deform slightly the projection as in figure 4, where vertical 
lines indicate complex points sent to same real point. That is,
if we deform near the imaginary axis, positive (resp. negative) 
imaginary numbers are sent to positive (resp. negative) real numbers.

\vskip 10pt

\ilustracion 20,40(trenza1){}.

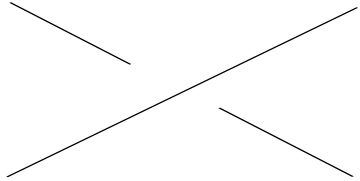
\ilustracion 20,40(trenza2) Braid corresponding to $\sqrt{z}$
around the unit circle.

\bigbreak
We end this section with some comments on
the relationship of the fundamental group
of affine and projective curves. Let
$\C\subset\bp^2$ be a reduced projective
curve of degree $d$. Let $G$ be the fundamental 
group $\pi_1(\bp^2\setminus\C)$.

Let $\Li$ be a line in $\bp^2$ transversal to $\C$. 
Let us identify $\bc^2$ with $\bp^2\setminus \Li$ and 
set $\Caf:=\C\setminus \Li$, which is an affine curve.
Let us denote by $H$ the fundamental group of the affine
curve, that is, $H=\pi_1(\bc^2\setminus\Caf)$. Let 
$i_*\:H\to G$ be the homomorphism induced by the inclusion. 
It is easily seen that this homomorphism is surjective.
With the notation of this section, the kernel of this mapping
is the cyclic infinite central subgroup $K$ of $H$
generated by $\mu_d\cdot\dots\cdot\mu_1,$ the meridian of the
line at infinity. These facts are easy consequences of the 
Zariski-Van Kampen method.

If one drops the assumption of transversality on $\Li$, the
product $\mu_d\cdot\dots\cdot\mu_1$ is no longer a meridian
of $\Li$. One can still describe a meridian of $\Li$ as a word
in $\mu_1,\dots,\mu_d$. Note that such meridian is not 
necessarily central in $H$. An example of this situation will
be carried out in \S 3. For a more general description see
\cite{A et al.}.

These groups are generated by the meridians of the curves.
Let $\mu\:G\to\bz/d\bz$ be the group homomorphism defined by
sending each meridian to $[1]$. Let us consider the quotient
$\pi\:\bz\to\bz/d\bz$ and the pull-back diagram for $\mu$
and $\pi$:
$$\matrix
\tilde H\ \ &\miya{{\tilde\mu}} &\bz\ \ \\
\downarrow\tilde\pi&    &\downarrow\pi \\
G\ \ & \miya{{\mu}} & \bz/d\bz.
\endmatrix$$

\newbox\genline
\global\setbox\genline=\etiqueta

\proclaim{Proposition\cita} There exists
an isomorphism $\varphi\:H\to\tilde H$ such
that $\tilde\pi\circ\varphi=i_*$. In particular,
the pair $(\bp^2,\C)$ determines $H$ providing
$\Li$ is transversal to $\C$.
\endproclaim

\demo{Proof} We know that 
$$\tilde H\:\{(\alpha,x)\in G\times\bz|\ 
\mu(\alpha)\equiv x\ \bold{mod}\ d\}.$$
In particular, it is the kernel of
$$\matrix
G\times\bz&\longrightarrow&\bz/d\bz\\
(\alpha,x)&\mapsto&\mu(\alpha)-[x].
\endmatrix$$
The remarks above concerning Zariski-Van Kampen method
and Reidemeister-Schreier method give the result. \qed
\enddemo

\finparrafo

\head\sec Computation of braid monodromy for $\C^{(1)}$
and $\pi_1(\PP^2 \setminus \C^{(1)})$
\endhead

Consider the equation of $\C^{(1)}$ --\,example \copy\ejemi\,--
in the affine chart $\{z \not = 0\}$. The affine curve will
be denoted by $\Caf^{(1)}$ and the discriminant set with
respect to $y$ by $S_1$. It is easily seen that $\Caf^{(1)}$
has only real non-transversal lines --\,figure 1\,-- and that
$\#S=6$. We order these points $x_1>\cdots>x_6$ and choose
$x_0>x_1$. In what follows we will use the path $\Gamma$ and 
the disks $D_i$ of the standard construction described in the
previous sections.

\medbreak
\ilustracion 18,120(base) Standard basis for $\bc$ and $S_1$.

From figure 1 we list in a table the braids associated with 
each path. We have followed the conventions:

\smallbreak\item{--} 
The thick lines correspond to the conic
$\C^{(1)}_1$. The fat lines correspond to the quartic
$\C^{(1)}_2$.

\smallbreak\item{--} 
The dotted lines correspond to the real 
part of non-real complex conjugate points of the curves with 
real first coordinate.

\smallbreak\item{--} 
At the intersections of these dotted lines
we put a thin string if the imaginary part of the points
corresponding to the conic have a bigger absolute value.
Otherwise, we put a thick string.

Let $\F:=\{x=x_0\}$, where $x_0\gg0$, be a generic vertical line.
Let also choose generators for the fundamental group of 
$\check \F:=\F\setminus \Caf^{(1)}$. From right to left we call 
them $a_1$, $a_2$, $a_3$ and $a_4$. Note that $a_1$ and $a_4$ are 
meridians of $\Cqaf^{(1)}$ whereas $a_2$ and $a_3$ are meridians 
of $\Ccaf^{(1)}$. Also note that $a_2$ corresponds to a point 
with positive imaginary part.

\vskip 10pt

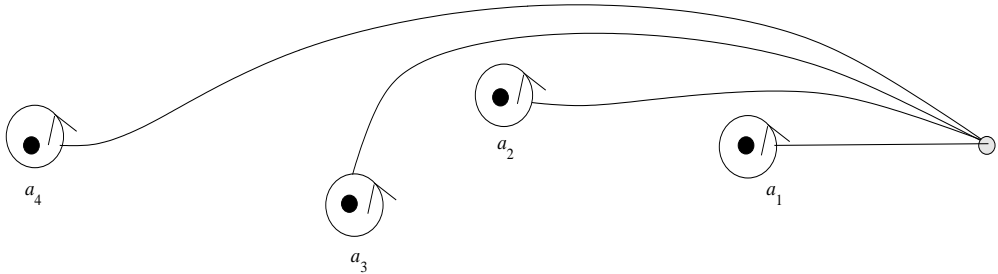
\ilustracion 30,110(gens) Generic fiber with generators of the
fundamental group homotopic to a standard basis.

The braids are:

$$\vbox{\offinterlineskip
\halign{\tv\centro{#}&\tv\centro{#}&\tv#\cr
\hline
Paths & Braids &\cr
\hline
$\Gamma_1$ & $\sigma_3^{-1}\sigma_2$ &\cr \hline
$\alpha_1$ & $\sigma_1$ &\cr \hline
$\alpha_1^+$ & $1$ &\cr \hline
$\Gamma_2$ & $\sigma_2^{-1}(\sigma_3^{-1}\sigma_1)^2\sigma_2$ &\cr \hline
$\alpha_2$ & $\sigma_3$ &\cr \hline
$\alpha_2^+$ & $\sigma_3$ &\cr \hline
$\Gamma_3$ & $\sigma_2^{-1}\sigma_1$ &\cr \hline
$\alpha_3$ & $\sigma_2^2$ &\cr \hline
$\alpha_3^+$ & $\sigma_2$ &\cr \hline
$\Gamma_4$ & $\sigma_1^{-1}\sigma_2$ &\cr \hline
$\alpha_4$ & $\sigma_1$ &\cr \hline
$\alpha_4^+$ & $\sigma_1$ &\cr \hline
$\Gamma_5$ & $1$ &\cr \hline
$\alpha_5$ & $\sigma_2^{16}$ &\cr \hline
$\alpha_5^+$ & $\sigma_2^8$ &\cr \hline
$\Gamma_6$ & $1$ &\cr \hline
$\alpha_6$ & $\sigma_3$ &\cr \hline
}}
$$

The braid monodromy is:
$$
\matrix
\Phi_1\:&\pi_1(\bc^2\setminus S_1)&\to&B_4\\
&\gamma_1&\mapsto&(\sigma_3^{-1}\sigma_2)\cdot\sigma_1\\
&\gamma_2&\mapsto&(\sigma_3^{-1}\sigma_2\sigma_2^{-1}
\sigma_3^{-2}\sigma_1^2\sigma_2)\cdot\sigma_3\\
&\gamma_3&\mapsto&(\sigma_3^{-1}\sigma_2\sigma_2^{-1}
\sigma_3^{-2}\sigma_1^2\sigma_2\sigma_3
\sigma_2^{-1}\sigma_1)\cdot\sigma_2^2\\
&\gamma_4&\mapsto&(\sigma_3^{-1}\sigma_2
\sigma_2^{-1}\sigma_3^{-2}\sigma_1^2\sigma_2\sigma_3
\sigma_2^{-1}\sigma_1\sigma_2\sigma_1^{-1}\sigma_2)\cdot\sigma_1\\
&\gamma_5&\mapsto&(\sigma_3^{-1}\sigma_2\sigma_2^{-1}
\sigma_3^{-2}\sigma_1^2\sigma_2\sigma_3
\sigma_2^{-1}\sigma_1\sigma_2\sigma_1^{-1}\sigma_2\sigma_1)
\cdot\sigma_2^{16}\\
&\gamma_6&\mapsto&(\sigma_3^{-1}\sigma_2\sigma_2^{-1}
\sigma_3^{-2}\sigma_1^2\sigma_2\sigma_3
\sigma_2^{-1}\sigma_1\sigma_2\sigma_1^{-1}\sigma_2\sigma_1\sigma_2^8)
\cdot\sigma_3
\endmatrix
$$
where $a\cdot b:=aba^{-1}$.

This braid monodromy allows us to construct a presentation of
the fundamental group of the complement of the affine curve.
Our goal is to compute $\pi_1(\bp^2\setminus \C^{(1)})$.
In order to do this it is enough to add any meridian of $\Li$
as a relation.

\medbreak
In order to construct a meridian of the line at infinity, we blow
up the projection point (the $\ba_3$ singular point).
The result is a ruled surface $\FF_1$ with a $(-1)$-section $E$, where
the strict transform $M$ of the line at infinity is one of the fibers
of the ruling (the other ones being the compactified vertical lines).
The real picture near $M$ is shown in the figure below.

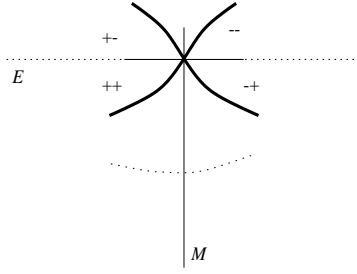
\ilustracion 30,40(infinito) Dotted lines are real parts of the two
branches of the conic.

Take the boundary of a tubular neighborhood $\T$ of $E$ in 
$\FF_1$. It is an $S^1$-bundle $\pi\:\T\to\bp^1$ over $\bp^1$ 
with Euler number equal to $-1$. This means the following:
take disks $\Delta_1,\Delta_2\subset\bp^1$ which cover $\bp^1$ and
whose intersection is the common boundary and take a section 
$S:\Delta_1\to \T$. Let $a\in\partial\Delta_1$ and $b:=S(a)$. 
Let $s$ be the path which runs counterclockwise 
$\partial S(\Delta_1)$ and is based at $b$. Let $e$ be the 
oriented fiber over $a$ based at $b$. Then, note that $s$ and $e$ 
are homotopic on the solid torus $\pi^{-1}(\Delta_2)$.

On the base $\PP^1$, we choose a (complex) coordinate $u$ 
centered at $\infty$ and on the fiber we choose a (complex) 
coordinate $v$ centered at $E$. Note that we can do this 
preserving real coordinates. We show the four quadrants 
of the real parts in figure 8. We can suppose that near 
$M\cap E$ (which is the origin of this system) $\T$ is the 
closed polydisk of radius $2\varepsilon$ for some 
$\varepsilon>0$. Fix the point $(\varepsilon,\varepsilon)$ as 
base point for the various fundamental groups. Note that
$\Delta_1$ can be chosen to be the disk of radius $\varepsilon$
in the $u$-coordinate, and $S$ the section 
$u\mapsto(u,\varepsilon)$. Preserving the notation introduced 
above, $s=e$ in the fundamental group of a solid torus. We note 
that this solid torus (contained in the original copy of $\bc^2$) 
does not intersect the curve $\Caf^{(1)}$ if $\T$ is chosen small 
enough. Then, one has $s=e$ in $\bc^2\setminus \Caf^{(1)}$. Let us 
look at the relative position of the image of $S$, $\Caf^{(1)}$ and 
$M$. Figure 9 shows the situation in the disk $v=\varepsilon$.

\vskip 10pt

\ilustracion 30,30(seccion) Situation in $S$.

The path $m$ is a meridian of $M$, and the paths $c_1$ and $c_2$ 
are meridians of $\C^{(1)}$. In $\bc^2\setminus \Caf^{(1)}$ one 
has the equality $e=s=c_2mc_1$. Note that $c_1$ is homotopic in
$\bc^2\setminus \Caf^{(1)}$ to a path in the line $\F_\varepsilon$ 
of equation $x=\varepsilon^{-1}$. We can suppose that 
$\F_\varepsilon=\F$, i.e., $\varepsilon x_0=1$. It is easily seen 
that $c_1=a_1$, where $a_1,\dots,a_4$ is the standard basis in $\F$
--\,see figure 7.

Let us choose a standard basis on the line $\F_{-\varepsilon}$ of 
equation $x=-\varepsilon^{-1}$. Let $\beta$ be the dotted path 
in figure 9. Let us take a standard basis in $\F_{-\varepsilon}$
and denote by $a_1',\dots,a_4'$ the paths obtained by applying the 
mapping $*\mapsto\beta\cdot*\cdot\beta^{-1}$ to this basis. It is 
easily seen that $c_2$ is homotopic to $a_1'$ in 
$\bc^2\setminus \Caf^{(1)}$.

There is a braid relating these two basis in $\bc^2\setminus \Caf^{(1)}$:
$$a_i'=(a_i)^{\sigma_1\sigma_2\sigma_3\sigma_2\sigma_1},\quad
i=1,\dots,4.$$
We deduce that $a_1'=a_4$. The properties of standard basis
imply that already in the fiber $\F=\F_\varepsilon$ we have
$a_4a_3a_2a_1e=1$. Then:
$$m=(a_1a_4a_3a_2a_1a_4)^{-1}.$$

Therefore a presentation for $\pi_1(\bp^2\setminus \C^{(1)})$
may be given as follows:

\smallbreak\item{--} Generators: $a_1,a_2,a_3,a_4$

\smallbreak\item{--} One class of relators of the type
$a_i=a_i^{\Phi(\gamma_j)}$, $i=1,\dots,4$, $j=1,\dots,6$.

\smallbreak\item{--} The relator $a_1a_4a_3a_2a_1a_4=1$.

In order to simplify this presentation, we have used the free 
software GAP (version 4, release 1, fix 6). In the appendix, 
we provide the GAP program we applied to get the result:

\newbox\gri
\global\setbox\gri=\etiqueta
\proclaim{Theorem\cita} The fundamental group $G_1$ of the complement
of $\C^{(1)}$ in $\bp^2$ has a presentation
$$|a,b:a^2(ab)^2=1,[a,b^2]=1|,$$
where $a$ is a meridian of the quartic and $b$ is a meridian of the 
conic. The abelianization exact sequence 
is a central extension of $\bz/2\bz$ by $\bz\times\bz/2\bz$.
Choosing $\omega$ the non-trivial element of $\bz/2\bz$ and $p,q\in G_1$ 
such that their images in $\bz\times\bz/2\bz$ generate this group,
the image of $q$ being of order $2$, then $q^2=pqp^{-1}q^{-1}=\omega$.
\endproclaim

\finparrafo

\head\sec Computation of braid monodromy for $\C^{(2)}$
and $\pi_1(\PP^2 \setminus \C^{(2)})$
\endhead

The procedure and notations will be according to the previous
section. Let $S_2$ be the discriminant set in this case. Note 
--\,figure 2\,-- that $S_2$ consists of four points. Let us 
choose the standard fiber $\F:=\{x=x_0\}$ for $x_0\gg 0$ and 
denote the standard generators again by $a_1$, $a_2$, $a_3$ and 
$a_4$. Both $a_1$ and $a_4$ are meridians of $\Cqaf^{(2)}$,
whereas $a_2$ and $a_3$ are meridians of $\Ccaf^{(2)}$. Note
that $a_2$ corresponds to a point with positive imaginary part.

We use again a standard notation for the paths in the base:

$$\vbox{\offinterlineskip
\halign{\tv\centro{#}&\tv\centro{#}&\tv#\cr
\hline
Paths & Braids &\cr
\hline
$\Gamma_1$ & $\sigma_1^{-1}\sigma_3^{-1}$ &\cr \hline
$\alpha_1$ & $\sigma_2$ &\cr \hline
$\alpha_1^+$ & $1$ &\cr \hline
$\Gamma_2$ & $1$ &\cr \hline
$\alpha_2$ & $(\sigma_1\sigma_3)^7
(\sigma_2\sigma_1\sigma_3\sigma_2)$ &\cr \hline
$\alpha_2^+$ & $(\sigma_1\sigma_3)^4
(\sigma_2\sigma_1\sigma_3\sigma_2)$ &\cr \hline
$\Gamma_3$ & $1$ &\cr \hline
$\alpha_3$ & $\sigma_2$ &\cr \hline
$\alpha_3^+$ & $1$ &\cr \hline
$\Gamma_4$ & $\sigma_1^{-1}\sigma_3^{-1}$ &\cr \hline
$\alpha_4$ & $\sigma_2^2$ &\cr \hline
}}
$$

The braid monodromy is:
$$
\matrix
\Phi_1\:&\pi_1(\bc^2\setminus S_2)&\to&B_4\\
&\gamma_1&\mapsto&(\sigma_1^{-1}\sigma_3^{-1})\cdot\sigma_2\\
&\gamma_2&\mapsto&(\sigma_1^{-1}\sigma_3^{-1})\cdot((\sigma_1\sigma_3)^7
(\sigma_2\sigma_1\sigma_3\sigma_2))\\
&\gamma_3&\mapsto&((\sigma_1\sigma_3)^3
(\sigma_2\sigma_1\sigma_3\sigma_2))\cdot\sigma_2\\
&\gamma_4&\mapsto&((\sigma_1\sigma_3)^2
(\sigma_2\sigma_1\sigma_3\sigma_2))\cdot\sigma_2^2
\endmatrix
$$
\vskip .25truecm

The construction of a meridian $m$ of the line at infinity
is exactly as in the previous section and we obtain
$$m=(a_1a_4a_3a_2a_1a_4)^{-1}.$$

Since we have chosen as base fiber $\F$ the line of equation $x=x_0$,
we must translate this element to the standard basis $a_1,\dots,a_4$
in $\F$. Since
$$a_i'=a_i^{(\tau^{-1})},$$
we deduce that
$$m=(a_1a_4a_3a_2a_1a_3a_2a_3^{-1})^{-1}.$$

We construct a presentation of $\pi_1(\bp^2\setminus \C^{(2)})$
and we simplify it with GAP:

\newbox\grii
\global\setbox\grii=\etiqueta
\proclaim{Theorem\cita} The fundamental group $G_2$ of the complement
of $\C^{(2)}$ in $\bp^2$ has a presentation
$$|a,b:b^2=(ab)^4|,$$
where $a$ is a meridian of the quartic and $b$ is a meridian 
of the conic. The abelianization exact sequence is an 
extension of $\FF_3$ (the free group in three generators $r,s,t$) 
by $\bz\times\bz/2\bz$. One can choose $p,q\in G_2$ such that 
their images in $\bz\times\bz/2\bz$ generate this group,
the image of $q$ being of order $2$, and:
 
\smallbreak\item{--}$q^2=r$, $pqp^{-1}q^{-1}=t$.

\smallbreak\item{--}$prp^{-1}=ts$, $psp^{-1}=t^{-1}$,
$ptp^{-1}=tr^{-1}$.

\smallbreak\item{--}$qrq^{-1}=r$, $qsq^{-1}=tr^{-1}$,
$qtq^{-1}=rs$.

In particular, the groups $G_1$ and $G_2$ are not isomorphic.
This provides a proof of Theorem \copy\zarpair.
\endproclaim

\finparrafo

\head\sec From special to generic braid monodromy
\endhead

In the previous sections we have computed some braid monodromies
where the projection point is in the projective curve. In fact
the projection point is a tacnode of the curve and the line
at infinity is tangent to this tacnode. We will apply this computation
to find a generic braid monodromy in both cases. The idea is to slide
the projection point somewhere close to the former one but not on
the curve. In the first part, we study the situation near the infinity;
in the second part we concentrate the curves near the transversal 
lines to the primitive projection.

\example{Step 1} We keep the affine curve, but we move the projection
point along the line at infinity.
\medbreak
\ilustracion 43,90(gener1).

As in the introduction, the dotted lines represent the real part of the
complex conjugated roots. When two couples of roots have the same real part
the non-dotted line represents the one having biggest absolute value for
the imaginary part.
\endexample

\example{Step 2} We take the dotted vertical line as the new line
at infinity and preserve the projection.
\medbreak
\ilustracion 43,90(gener2).
\endexample

\example{Step 3} Finally we move the projection in
order to have a generic projection at the tacnode.
\medbreak
\ilustracion 43,90(gener3).
\endexample

We fix the dotted line as the generic fiber and
compute the braid monodromy for the part which
is to the right of the generic fiber. The paths 
are decomposed as in the previous sections:

$$\vbox{\offinterlineskip
\halign{\tv\centro{#}&\tv\centro{#}&\tv#\cr
\hline
Paths & Braids &\cr
\hline
$\tilde\Gamma_1$ & $1$ &\cr \hline
$\tilde\alpha_1$ & $\sigma_1$ &\cr \hline
$\tilde\alpha_1^+$ & $1$ &\cr \hline
$\tilde\Gamma_2$ & $\sigma_2^{-1}(\sigma_3^{-1}\sigma_1)\sigma_2$ 
&\cr \hline 
$\tilde\alpha_2$ & $\sigma_3$ &\cr \hline
$\tilde\alpha_2^+$ & $\sigma_3$ &\cr \hline
$\tilde\Gamma_3$ & $1$ &\cr \hline
$\tilde\alpha_3$ & $\sigma_4^4$ &\cr \hline
$\tilde\alpha_3^+$ & $\sigma_4^2$ &\cr \hline
$\tilde\Gamma_4$ & $1$ &\cr \hline
$\tilde\alpha_4$ & $\sigma_5$ &\cr \hline
$\tilde\alpha_4^+$ & $1$ &\cr \hline
$\tilde\Gamma_5$ & $\sigma_4^{-1}\sigma_5\sigma_3^{-1}\sigma_4%
\sigma_2^{-1}\sigma_3\sigma_1^{-1}\sigma_2$ &\cr \hline
$\tilde\alpha_5$ & $\sigma_1$ &\cr \hline
}}$$

This part is common for both curves. Note that the projection 
of $\Caf^{(1)}$, coincides with the non-generic case. One can 
recover the generic monodromy just by shifting
$\sigma_i\mapsto\sigma_{i+1}$, for $i=1,2,3$, since a new thread
has appeared in the upper part.

For the curve $\Caf^{(2)}$, we may have the same braid monodromy
near the infinity and outside the point $\ba_{15}$.
For this point, the real deformation looks as follows.

\ilustracion 40,40(a15).

We recall that in the special monodromy, near the point $\ba_{15}$
all solutions were imaginary complex conjugated.  Since the projection
onto the real axis was not generic, we deformed this projection
in order such that the order of the points was
$((P')^+,P_4^+,P_4^-,(P')^-)$, where sub-index indicates de degree of 
the curve and the super-index indicates the sign of the imaginary part.
Let us consider the line in the special monodromy where the braid
$\Gamma_2$ starts. Take a line in this generic monodromy close to that 
line in order to start a braid $\Gamma_{2a}$. This braid will be very 
close to the braid $\Gamma_2$; then, one can only have interchanging of 
positive or negative points and finally a crossing to be in the situation
of last figure. Since there are only conjugated interchanges of positive
and negative points, one has the following:

$$\vbox{\offinterlineskip
\halign{\tv\centro{#}&\tv\centro{#}&\tv#\cr
\hline
Paths & Braids &\cr
\hline
$\Gamma_{2a}$ & $\sigma_1^{-2k}\sigma_3^{2k}\sigma_1\sigma_2$ &\cr \hline
$\alpha_{2a}$ & $\sigma_1$ &\cr \hline
$\alpha_{2a}^+$ & $\sigma_1$ &\cr \hline
$\Gamma_{2b}$ & $\sigma_2^{-1}\sigma_3$ &\cr \hline
$\alpha_{2b}$ & $\sigma_2$ &\cr \hline
$\alpha_{2b}^+$ & $\sigma_2$ &\cr \hline
$\Gamma_{2c}$ & $1$ &\cr \hline
$\alpha_{2c}$ & $\sigma_3^{16}$ &\cr \hline
$\alpha_{2c}^+$ & $\sigma_3^8$ &\cr \hline
}}
$$
Finally, since we performed a deformation, the following 
equality must be satisfied in the braid group:
$$\alpha_2^+=\Gamma_{2a}\alpha_{2a}^+
\Gamma_{2b}\alpha_{2b}^+\Gamma_{2c}\alpha_{2c}^+.$$
A straightforward computation gives:
$$\Gamma_{2a}\alpha_{2a}^+\Gamma_{2b}\alpha_{2b}^+
\Gamma_{2c}\alpha_{2c}^+=\sigma_1^{8-2k}\sigma_3^{2k}
\sigma_2\sigma_1\sigma_3\sigma_2.$$
Then $k=2$ is the only possible solution.
Replacing $\Gamma_2,\alpha_2,\alpha_2^+$ by these braids
and applying the shifting $\sigma_i\mapsto\sigma_{i+1}$
one obtains the generic braid monodromy for the affine
part.

\bigbreak
Note that, using the generic braid monodromy and 
Zariski-van Kampen method, one would be able to obtain a
presentation of $G_i$ in a more direct way, avoiding the 
calculations at the end of \S 3.

\finparrafo

\head\sec Characteristic varieties
\endhead

Notations and definitions in this section will mainly follow
\cite{L4} and \cite{Co}. We first want to recall what the 
characteristic varieties of a plane curve are. Let 
$\C=\C_1 \cup ... \cup \C_r$ be an algebraic curve in $\PP^2$ 
and consider $\Li$ a transversal line to $\C$. 
It is well known that $X=\PP^2 \setminus (\C \cup \Li)$ 
has the homotopy type of a finite 2-dimensional CW-complex and 
$H_1(X)=\ZZ^r$. We will denote by $\tilde{X}$ the universal abelian cover
of $X$. This has again the homotopy type of
a 2-dimensional CW-complex, but no longer finite. 
In order to study its homology, we can use the action of $H_1(X)$ over 
$\tilde{X}$, which makes any $H_k(\tilde{X})$ into a $\bz[H_1(X)]$-module.
Tensoring by $\bc$ we are led to study the structure of the modules
$H_k(\tilde{X},\bc)$ over the noetherian ring $\Lambda=\bc[\bz^r]$. We
define the $k$-th characteristic variety of $\C$ as follows
$$Char_k(\C):= Supp_{\Lambda}(\wedge^k H_1(\tilde{X},\bc)) =
Supp_{\Lambda}(\Lambda/F_k) \subset Spec \Lambda = (\bc^*)^r,$$
where $F_k$ is the $k$-th Fitting ideal of the module $H_1(\tilde{X},\bc)$.
In other words, $Char_k(\C)$ is a subvariety of the $r$-dimensional torus 
given by the zeroes of the $k$-th Fitting ideal of $H_1(\tilde{X},\bc)$.
Observe that $H_1(\tilde{X})=\pi_1(X)'/\pi_1(X)''$, and it is easy to
check that the aforementioned action of $H_1(X)=\pi_1(X)/\pi_1(X)'$ 
on $H_1(\tilde{X})$ corresponds to conjugation. As a result of this remark
we have that the varieties $Char_k(\C)$ are invariants of the fundamental 
group $\pi_1(X)$. In the case where $r=1$, the first characteristic 
varieties are nothing but the zeroes of the Alexander polynomial. 
In the non-irreducible case, the Fitting ideals $F_k$ are no longer 
principal and as a result
$Char_k(\C)$ will in general not be a hypersurface.

\smallbreak
Characteristic varieties can also be thought of in connection 
with homology with coefficients in rank one local systems of 
the fundamental group $\pi_1(X)$. This connection makes it 
possible to calculate the characteristic varieties of a curve 
$\C$ using the theory of adjoint and quasiadjoint ideals.

\smallbreak
Our purpose is to calculate $Char_k(\C^{(1)})$ and 
$Char_k(\C^{(2)})$. They will turn out to be different 
sets of points. Since the characteristic varieties are
invariants of the fundamental groups and by \copy\genline\ 
Theorem \copy\zarpair\ will be proved.

\smallbreak
One possible way to calculate the characteristic varieties is to 
compute the Fitting ideals $F_k$ by means of Fox calculus. For this
purpose one needs a finite presentation of the fundamental group. We 
will denote the Alexander modules of $\C^{(1)}$ and $\C^{(2)}$ by
$M^{(1)}$ and $M^{(2)}$ respectively. We will also consider the 
2-dimensional torus $Spec \Lambda = (\bc^*)^2$, where the first 
coordinate will represent a meridian around the quartic and the second 
coordinate will represent a meridian around the conic. 
Applying Fox calculus to the presentations 
given in Theorems \copy\gri\ and \copy\grii\ one can compute that
$$F_1(M^{(1)})=(t_1-1,t_2-1)$$
$$F_1(M^{(2)})=(t_1-1,t_2-1) \cdot
((t_2t_1+1)(t_2^2t_1^2+1),(t_2+t_2^2t_1+t_2^3t_1^2-1)).$$ 
Hence
$$Char_1(\C^{(1)})=\{(1,1)\}$$
$$Char_1(\C^{(2)})=\{(1,1),(-1,1),(\sqrt{-1},-1),(-\sqrt{-1},-1)\}.$$

\noindent
Another way to calculate characteristic varieties is by computing the
cohomology of the universal abelian covers $\tilde{X}^{(1)}$ and 
$\tilde{X}^{(2)}$ with coefficients in rank one local systems as 
mentioned before. These computations involve both local and global 
analyses and do not need a presentation of the fundamental group:

\smallbreak\item{$1.$} Local calculations

\smallbreak\itemitem{$a.$} Calculate the ideals of quasiadjunction 
for each singular point, the local polytopes of quasiadjunction and 
their faces.

\smallbreak\itemitem{$b.$} Put all the local data together in the 
global polytope of quasiadjunction.

\smallbreak\itemitem{$c.$} Determine the faces of the global 
polytope of quasiadjunction that are contributing faces. Each 
contributing face has associated with it an
integer number $l(\delta)$ called the level of $\delta$.

\smallbreak\item{$2.$} Global calculations

\smallbreak\itemitem{$a.$} Each contributing face $\delta$ has an 
ideal sheaf ${\Cal A}_{\delta}$ over $\PP^2$ associated with it. 
The contributing face will determine a subvariety of the $k$-th 
characteristic variety if and 
only if the irregularity of ${\Cal A}_{\delta}(d-3-l(\delta))$ is 
greater or equal to $k$, where $d$ is the degree of the curve $\C$ and
$l(\delta)$ is the level of the contributing face $\delta$. 

\smallbreak\itemitem{$b.$} Let $\delta$ be a contributing face from 
the previous step for
which the irregularity is positive, that is 
$\dim H^1({\Cal A}_{\delta}(d-3-l(\delta))) = k > 0$. 
Let $L_s(x_1,...,x_r)=\beta_s$ be a system of equations defining
$\delta \cup \overline{\delta}$ where 
$$\overline{\delta}=\{ \overline{1} - \overline{x} \mid 
\overline{x} \in \delta \}.$$
Suppose that each $L_s(x_1,...,x_r)$ is a linear form with integer
coefficients and that the g.c.d. of the maximal order minors in 
the matrix of coefficients is one. Then 
$$exp(2 \pi \sqrt{-1} L_s) = exp(2 \pi \sqrt{-1} \beta_s)$$
is the equation of a component of $Char_k(\C)$, where the coordinates 
in the $r$-dimensional torus $Spec \Lambda$ are given by
$t_i=exp(2 \pi \sqrt{-1} x_i)$. Moreover,
all non-trivial components of $Char_k(\C)$ can be obtained in this way.

In the following, we will work out each step for our curves $\C^{(1)}$ 
and $\C^{(2)}$.

\smallbreak\item{$1.$} The local computations can be carried out for both 
curves together, since they both have the same local configuration
of singularities.

\smallbreak\itemitem{$a.$} The ideals of quasiadjunction of a 
singular point $\ba_{15}$ given by two different components are 
$$m_{1,7} \subset m_{1,6} \subset ... \subset m_{1,2} \subset 
m_{1,1} = m \subset {\Cal O},$$
where $m_{1,i}=(y,x^i)$ ($i=1,...,7$), $y$ denotes the tangent line 
at the singular point and $x$ is independent to $y$. The local face
of quasiadjunction of $m_{1,i}$ is given by the equation 
$$8 x_1+8 x_2 = 7 - i$$
and hence the local polytope of quasiadjunction is

\vskip 10pt

\ilustracion 40,40(topo1).

The only proper ideal of quasiadjunction for a tacnode whose 
local branches belong to the same global irreducible component 
is the maximal ideal. There is only one face of quasiadjunction 
$x_1=1/4$ and it is associated with the maximal ideal.

Nodes have no ideals of quasiadjunction attached to them.

\smallbreak\itemitem{$b.$} The global polytope of quasiadjunction 
looks as follows, 

\vskip 10pt
\ilustracion 40,40(topo2).

\smallbreak\itemitem{$c.$} The contributing faces are the points 
in the intersection of faces in the above picture that belong to 
hyperplanes of type
$$4 x_1 + 2 x_2 = l \in \bn.$$
Hence the only contributing face is
$\delta = \{(\frac{1}{4},\frac{1}{2})\}$ and its level $l(\delta)=2$.

\smallbreak\item{$2.$} Global calculations do depend on the curve. 
We will distinguish each case by superindices.

\smallbreak\itemitem{$a.$} The ideal sheaf associated with $\delta$ 
can be defined as follows 
$$({\Cal A}^{(i)}_{\delta})_P=
\cases
m_{1,2} & \text{$P$  is the singularity $\ba_{15}$
in $\C^{(i)}$}\\
m & \text{if $P$ is the singularity $\ba_{3}$ 
in $\C^{(i)}$}\\
{\Cal O}_{\PP^2,P} &\text{otherwise}
\endcases
$$
Let us calculate the irregularity of both ideal sheaves.
$$\dim H^1({\Cal A}^{(1)}_{\delta}(6-3-2))=
\dim H^0({\Cal A}^{(1)}_{\delta}(1)) - 
\chi({\Cal A}^{(1)}_{\delta}(1)) = 0$$
since 
$$\chi({\Cal A}^{(1)}_{\delta}(1)) = \chi({\Cal O}_{\PP^2}(1)) 
- \dim(m_{1,2} \oplus m)= 3 - 3 = 0.$$
Analogously
$$\dim H^1({\Cal A}^{(2)}_{\delta}(6-3-2))=
\dim H^0({\Cal A}^{(2)}_{\delta}(1)) - 
\chi({\Cal A}^{(2)}_{\delta}(1)) = 1.$$

\smallbreak\itemitem{$b.$} Hence $Char_k(\C^{(1)})$ has at most 
only trivial components, whereas $Char_1(\C^{(2)})$ has a 
non-trivial component coming from the face
$\delta$. The defining equations of $\delta$ are
$$L_1(x_1,x_2) = 2 x_2 = 1 = \beta_1$$
and
$$L_2(x_1,x_2) = 2 x_1 + x_2 = 1 = \beta_2.$$
Hence,
$$Char_1(\C^{(2)})=\{ t_2^2=1, t_1^2t_2=1 \}=$$
$$=\{(1,1),(-1,1),(\sqrt{-1},-1),(-\sqrt{-1},-1)\}$$

These computations provide another proof for \copy\zarpair.
\finparrafo

\head\sec Dihedral covers
\endhead

Let $Y$ be a smooth projective variety and let $X$ be a normal 
variety with finite morphism $\pi : X \to Y$. We denote the 
rational function fields of $X$ and $Y$ by $\CC(X)$ and $\CC(Y)$
respectively. The field $\CC(X)$ is a finite extension of $\CC(Y)$ 
with $[\CC(X):\CC(Y)] = \deg\pi$. We call $X$ a $\D_{2n}$ cover if 

\smallbreak\item{(i)} 
$\CC(X)$ is a Galois extension of $\CC(Y)$ and

\smallbreak\item{(ii)}  its Galois
group is the dihedral group 
$\D_{2n} = \langle \sigma, \tau | \ 
\sigma^2 = \tau^n = (\sigma\tau)^2 =1 \rangle$. 

\smallbreak
Note that, given a $\D_{2n}$ cover $X$, a double cover of $Y$ is 
canonically determined by considering the $\tau$-invariant
field of $\CC(X)$. We denote this double cover by $D(X/Y)$.
In \cite{T2}, the last author studied such covers, and proved
the following result for $\D_{2n}$ ($n$ even) covers.

\medbreak

\newbox\dihedral
\global\setbox\dihedral=\etiqueta

\proclaim{Proposition\cita} 
Let $Y$ be a simply connected 
smooth projective variety and let $f :Z \to Y$ be a smooth 
double cover, and $\sigma$ denotes the covering transformation. 
Suppose that there exist effective divisors,
$D_1$ and $D_2$ on $Z$ satisfying the following conditions:

\smallbreak\item{\rm(i)} $D_1$ and $\sigma^*D_1$ have no common 
components; and if we let $D_1 = \sum_i a_i D_{i, 1}$ denote the 
irreducible decomposition of $D_1$, then the greatest common 
divisor of $a_i$'s and $n$ is one.

\smallbreak\item{\rm(ii)} If $D_2$ is not empty, $D_2$ is a 
reduced divisor and there exists a divisor, $B_2$, on $Y$ 
such that $D_2 = f^*B_2$.

\smallbreak\item{\rm(iii)} There exists a line bundle $\cL$ 
on $Z$ such that $D_1 + \frac n2 D_2 - \sigma^*D_1 \sim n \cL$ 
($n$ even).

Then there exists a $\D_{2n}$ cover, $X$, of $Y$ such that

\smallbreak\item{\rm(a)} $D(X/Y) Z$ and 

\smallbreak\item{\rm(b)} $\Delta(X/Z) \subset
\Supp (D_1 + \sigma^*D_1 + D_2)$.
\endproclaim

\smallbreak

\demo{Proof} By Remark 3.1, in \cite{T2} and a similar 
argument to the proof of Proposition 1.1 in \cite{T3}, 
Proposition 1.6 follows.\qed
\enddemo

The converse of Proposition \copy\dihedral\ 
in the following sense also holds.

\smallbreak

\newbox\dihedralconv
\global\setbox\dihedralconv=\etiqueta

\proclaim{Proposition\cita} 
Let $\pi : X \to Y$ be a $\D_{2n}$ 
cover such that $D(X/Y)$ is smooth. Then there exist (possibly empty) 
effective divisors, $D_1$ and $D_2$, and a line bundle $\cL$ on $D(X/Y)$ 
satisfying the following conditions:

$(i)$ If $D_1 \neq \emptyset$, $D_1$ and $\sigma^*D_1$ have 
no common components.

$(ii)$ If $D_2$ is not empty, $D_2$ is a reduced divisor and there 
exists a divisor, $B_2$, on $Y$ such that $D_2 = f^* B_2$.

$(iii)$ $D_1 + \frac n2 D_2 - \sigma^*D_1 \sim n \cL$.

$(iv)$ $\Delta(X/D(X/Y)) = \Supp(D_1 + \sigma^*D_1 + D_2)$.
\endproclaim

See \cite{T2} for a proof.

\newbox\dneven
\global\setbox\dneven=\etiqueta

\subhead\ssec $\D_{2n}$ ($n$ even) covers of $\PP^2$ \ 
\endsubhead

\smallbreak

In this subsection, we confine ourselves to studying a $\D_{2n}$ 
($n$: even) cover, $S$, of $\PP^2$ satisfying the following properties:

$(\ast)$ $\Delta(S/\PP^2)$ consists of a reduced curve, $B$, of
even degree and a line, $\Li$, such that

$(i)$ $\Delta(D(S/\PP^2))=B$,

$(ii)$ $B$ has at most simple singularities, and

$(iii)$ $B$ meets $\Li$ transversally.

\smallbreak

Let $f' : Z' \to \PP^2$ be the double cover branched at $B$,
where $f'$ is finite and $Z'$ is normal. Let
$\mu\:Z\to Z'$ be

$$\matrix
Z' & \hiya{{\mu}} & Z \\
f' \downarrow &   & \downarrow \\
\PP^2 & \hiya{q} & \Sigma
\endmatrix
$$
the canonical resolution, where $q : \Sigma \to \PP^2$ is a
sequence of blowing-ups so that the induced morphism
$Z \to \Sigma$ gives a smooth finite double cover. We denote
$f' \circ \mu$ by $f$. Let $\NS(Z)$ be the N\'eron-Severi group of $Z$.
Let $T$ be the subgroup of $\NS(Z)$ generated by $f^*\Li$ and all the
irreducible components of the exceptional divisor of $\mu : Z \to Z'$.

Under the assumption $(\ast)$, $Z$ is simply connected. Hence
$\NS(Z) = \Pic(Z)$ and $\NS(Z)$ is a lattice with respect to
the intersection pairing. Also $T$ is regarded as a sublattice
of $\NS(Z)$ and has a orthogonal decomposition
$$
T = \ZZ f^*\Li \bigoplus \bigoplus_{x \in \Sing(Z')}
\left (\bigoplus_i \ZZ \Theta_{i, x} \right ),
$$
where $\Theta_{i,x}$ denotes an irreducible component of the
exceptional divisor arising from $x$. We
recall that $\Sing(Z')$ is in bijection with $\Sing(B)$.

Under this notation, we have

\smallbreak

\newbox\ntorsion
\global\setbox\ntorsion=\etiqueta

\proclaim{Proposition\cita}
If there exits a $\D_{2n}$ cover, 
$S$, of $\PP^2$ such that 

\smallbreak\item{\rm(i)} $\Delta(S/\PP^2)=B+\Li$ and 

\smallbreak\item{\rm(ii)} 
$D(S/\PP^2) = Z'$.

Then $\NS(Z)/T$ has $n$-torsion.
\endproclaim

{\it Proof.} Let $\tS$ be the $\CC(S)$-normalization of 
$\Sigma$ and let $\tpi : \tS \to \Sigma$ be the induced morphism. 
Then $\tpi : \tS \to \Sigma$ satisfies:

$(i)$ $\tpi : \tS \to \Sigma$ is a $\D_{2n}$ cover of $\Sigma$,

$(ii)$ $D(\tS/\Sigma) = Z$, and

$(iii)$ the branch locus of $\tS \to Z$ is contained in $f^*\Li$ 
and the exceptional set of $\mu$.

Let $D_1$ and $D_2$ be the effective divisors on $Z$ as in 
Proposition\copy\dihedralconv. Then by the proof of Proposition 0.7, 
\cite{T2}, there exists a rational function
$\varphi$ such that

(a) $(\varphi) = D_1 + \frac n2 D_2 - \sigma^*D_1 + nD_0$, where
$D_0$ is some divisor on $Z$, and

(b) the polynomial $X^n - \varphi$ is irreducible in $\CC(Z)[X]$, 
i.e. $\CC(S)=\CC(Z)(\sqrt[n]{\varphi})$.

Since $\Supp (D_1 + \sigma^*D_1 + D_2) = \Delta(S/Z)$ and $(iii)$, 
the equality in (a) shows that $D_0$ gives rise to an element, 
$\alpha$, of $\NS(S)/T$ whose order divides $n$. Our assertion 
follows from the following claim:

\smallbreak

\proclaim{Claim} 
The order of $\alpha$ is $n$.
\endproclaim

\demo{Proof of Claim} Suppose that the order of $\alpha$ is 
$d$ and put $d_1 = n/d$. Hence $dD_0 \in T$ and therefore 
$dD_0 \sim af^*\Li + \sum_x\sum_i b_{i, x} \Theta_{i, x}$. 
As we may assume that
$$\matrix
D_1 & = & \sum_{i, x}a_{i, x} \Theta_{i, x}, \quad 0 \le a_{i, x} < n, \\
D_2 & = & m f^*\Li + \sum_{i,x}m_{i, x}\Theta_{i, x}, \quad 
m, m_{i,x} \text{\ is \ either $0$ or $1$.},
\endmatrix$$
we have
$$
\multline
D_1 + \frac n2 D_2 - \sigma^* D_1 + n D_0 =
\sum_{i, x}a_{i, x} \Theta_{i, x} + \frac n2 (m f^*\Li +
\sum_{i,x}m_{i,x} \Theta_{i, x})-\\
 - \sum_{i, x}a_{i, x} 
\sigma^*\Theta_{i, x} + d_1 af^*\Li + 
d_1\sum_{i, x}b_{i,x}\Theta_{i,x} \sim 0.
\endmultline
$$
Since $D_1$, $\sigma^*D_1$ and $D_2$ have no common components, the 
above relation is non-trivial, unless all the coefficients of $f^*\Li$ and
$\Theta_{i,x}$'s are $0$. This implies that $d_1$ divides $n/2$ and all
$a_{i, x}$'s, and we can write $D_1=d_1D_1'$ for some effective divisor
$D'_1$ on $Z$. Hence we have
$$d_1(D_1' + \frac n{2d_1} D_2 - 
\sigma^*D_1' + \frac n{d_1}D_0) \sim 0.$$
As $Z$ is simply connected, this is reduced to
$$D_1' + \frac n{2d_1} D_2 - \sigma^*D_1' - \frac n{d_1}D_0 \sim 0.$$
Let $\psi$ be a rational function on $Z$ such that
$$(\psi)=D_1' + \frac n{2d_1}D_2 - \sigma^*D_1' + \frac n{d_1}D_0.$$
Thus $\varphi/\psi^{d_1}$ is a constant. Hence, multiplying 
by a suitable constant, we have $\varphi = \psi^{d_1}$. But this 
contradicts the condition (b) about $\varphi$.\qed
\enddemo

This finishes the proof of Proposition \copy\ntorsion. \qed

\finsparrafo

\subhead\ssec Non-existence of a certain $\D_{16}$ cover of 
$\PP^2$ branched at $\C^{(1)}+\Li$ \ 
\endsubhead
\smallbreak

We keep the notation from \S 1.
Let $\C^{(1)}$ be the sextic from example \copy\ejemi.
In this section, we disprove the existence of a $\D_{16}$ cover, 
$S$, of $\PP^2$ such that

(i) $\Delta(S/\PP^2) = \C^{(1)} + \Li$ and

(ii) $D(S/\PP^2) = Z'$.

Suppose that such a cover exists. Then  $\NS(Z)/T$ has eight
torsion and $\NS(Z)/T$ is a torsion group. Denoting the 
absolute value of the determinant of the intersection 
matrix of $\NS(Z)$ (resp. $T$) by 
$\disc (\NS(Z))$ (resp. $\disc T$), we have
$$\disc(\NS(Z)) \le \frac {\disc (T)}{64}=4.$$
On the other hand, we have the following lemma which leads us to a
contradiction; and no $\D_{16}$ cover as above exists.

\smallbreak

\proclaim{Lemma\cita} 
$\disc (\NS(Z))=16.$
\endproclaim

\demo{Proof} Let $X_{8211}$ be the rational elliptic surface 
in \copy\ellipt, and let $\hat {X}_{8211}  \to X_{8211}$ be 
blowing-ups at two nodes of 2 $I_1$ fibers. By what we have 
seen in \copy\dneven, $\Sigma$ coincides with $\hat {X}_{8211}$, 
and $Z$ is a double cover branched at the proper transforms 
of $I_1$ fibers. Hence there exists an elliptic fibration on 
$\varphi : Z \to \PP^1$ with a section having singular fibers 
$4I_2$, $2I_8$. This implies that $\disc (\NS(Z))=16$ by 
\cite{Mi-P2}.\qed
\enddemo

\finsparrafo

\subhead\ssec Existence of a certain $\D_{16}$ cover 
branched at $\C^{(2)} + \Li$ \ 
\endsubhead
\smallbreak

Let $\C^{(2)}$ be the sextic curve described in \copy\ejemii, 
and let $\L_1$ be the tangent line at $\ba_{15}$. The strict 
transform of $\L_1$ splits up into two components, $L_1^+$ and 
$L_1^-$. We label the irreducible components of the exceptional 
divisor of $\mu$ as follows:

\vskip 20pt

\ilustracion 50,80(toku).

With this notation, we have

\smallbreak

\newbox\lplus
\global\setbox\lplus=\etiqueta

\proclaim{Lemma\cita} 
$$
L_1^+ \sim_{\QQ} \frac 12 f^*\Li - \frac 18 \left ( 7\Theta_{1,1} +
\sum_{k=2}^{15}(16 - k)\Theta_{k, 1} \right ) - \frac 14 
\sum_{k=1}^3(4 - k) \Theta_{k, 2}.
$$
\endproclaim

\demo{Proof} Put
$$
D_{\alpha} = L_1^+ - \frac 12 f^*\Li + \frac 18 \left ( 7\Theta_{1,1} +
\sum_{k=2}^{15}(16 - k)\Theta_{k, 1} \right ) + \frac 14 
\sum_{k=1}^3(4 - k) \Theta_{k, 2}.
$$
By \cite{T4}, Lemma 1.2 and Lemma 1.4, $D_{\alpha}$ satisfies that

(i) $D_{\alpha}=0 \ \mod \ T \otimes \QQ$ and

(ii) $D_{\alpha} \perp T$ with respect to the intersection pairing.

We can easily check $D_{\alpha}^2=0$. Hence $D_{\alpha} \approx_{\QQ} 0$
by Lemma 1.3 \cite{T4}. Since $Z$ is simply connected, 
we are done.\qed
\enddemo

\newbox\lequiv
\global\setbox\lequiv=\etiqueta

\proclaim{Lemma\cita} 
Put
$$\matrix
D_1 & = & \sum_{k=1}^7(8 - k)\Theta_{k, 1} + 6\Theta_{1,2}, \\
D_2 & = & f^*\Li + \Theta_{2,2} \\
D_0 & = & L_1^+ - f^*\Li + \Theta_{15,1} + \Theta_{8,1} +
\sum_{k=1}^7(\Theta_{k,1} + \Theta_{16-k, 1}) + \Theta_{3,2}
\endmatrix$$
Then $D_1$, $D_2$ and $D_0$ satisfies
$$
D_1 + 4D_2 - \sigma^*D_1 + 8D_0 \sim 0.
$$
\endproclaim

\demo{Proof} By straightforward computation and Lemma \copy\lplus, 
we have Lemma \copy\lequiv. \qed
\enddemo

\medbreak

The three divisors in Lemma \copy\lequiv\ satisfy the following 
conditions:

(i) $D_1$ and $D_2$ are effective divisors satisfying the conditions 
in Proposition \copy\dihedral.

(ii) $\Supp (D_1 + \sigma^*D_1 + D_2)$ is contained in
$\Supp(f^*\Li)$ and the exceptional set of $\mu$.
Hence by Proposition \copy\dihedral, we infer that there exists a 
$\D_{16}$ cover, $S$,  of $\PP^2$ branched at $\C^{(2)}+\Li$ 
and having $D(S/\PP^2)=Z'$.

\medbreak

From the results obtained in this section, 
one has Theorem \copy\zarpair.

\head Appendix: GAP programs
\endhead

Both programs have common beginnings and ends:
\bigbreak
\begingroup
\tt
\parindent=0cm
\obeylines
g:=FreeGroup(4,"a");
lista:=GeneratorsOfGroup(g);
a1:=lista[1];
a2:=lista[2];
a3:=lista[3];
a4:=lista[4];
lista1:=[a2,a2*a1/a2,a3,a4];
k1:=GroupHomomorphismByImages(g,g,lista,lista1);;
lista1m:=[a2\^{ }a1,a1,a3,a4];;
k1m:=GroupHomomorphismByImages(g,g,lista,lista1m);;
lista2:=[a1,a3,a3*a2/a3,a4];;
k2:=GroupHomomorphismByImages(g,g,lista,lista2);;
lista2m:=[a1,a3\^{ }a2,a2,a4];;
k2m:=GroupHomomorphismByImages(g,g,lista,lista2m);;
lista3:=[a1,a2,a4,a4*a3/a4];;
k3:=GroupHomomorphismByImages(g,g,lista,lista3);;
lista3m:=[a1,a2,a4\^{ }a3,a3];;
k3m:=GroupHomomorphismByImages(g,g,lista,lista3m);;
id:=GroupHomomorphismByImages(g,g,lista,lista);;

igual:=function(elemento,trenza) 
\qquad return Image(trenza,elemento)/elemento; 
end;
\bigbreak
\centerline{($\cdots$)}
\bigbreak
h:=g/rel;
P:=PresentationFpGroup(h);
TzGoGo(P);
TzPrintRelators(P);
\par
\endgroup

\bigbreak
The central part for $G_1$ is:

\begingroup
\tt
\parindent=0cm
\obeylines
\bigbreak
\centerline{($\cdots$)}
\bigbreak

camino:=k2m*k1m\^{ }2*k3\^{ }2*k2*k2m*k3;
caminom:=k3m*k2*k2m*k3m\^{ }2*k1\^{ }2*k2;

camino1:=camino*k3m*k2;;
camino1m:=k2m*k3*caminom;;
sing1:=k1;;
trenza1:=camino1*sing1*camino1m;;
rel:=List(lista,u->igual(u,trenza1));;

camino2:=camino1*k2m*k3m\^{ }2*k1\^{ }2*k2;;
camino2m:=k2m*k1m\^{ }2*k3\^{ }2*k2*camino1m;;
sing2:=k3;;
trenza2:=camino2*sing2*camino2m;;
rel:=Union(rel,List(lista,u->igual(u,trenza2)));;

camino3:=camino2*k3*k2m*k1;;
camino3m:=k1m*k2*k3m*camino2m;;
sing3:=k2\^{ }2;;
trenza3:=camino3*sing3*camino3m;;
rel:=Union(rel,List(lista,u->igual(u,trenza3)));;

camino4:=camino3*k2*k1m*k2;;
camino4m:=k2m*k1*k2m*camino3m;;
sing4:=k1;;
trenza4:=camino4*sing4*camino4m;;
rel:=Union(rel,List(lista,u->igual(u,trenza4)));;

camino5:=camino4*k1;;
camino5m:=k1m*camino4m;;
sing5:=k2\^{ }16;;
trenza5:=camino5*sing5*camino5m;;
rel:=Union(rel,List(lista,u->igual(u,trenza5)));;

camino6:=camino5*k2\^{ }8;;
camino6m:=k2m\^{ }8*camino5m;;
sing6:=k3;;
trenza6:=camino6*sing6*camino6m;;
rel:=Union(rel,List(lista,u->igual(u,trenza6)));;

rel:=Union(rel,[Image(caminom,a1*a4*a3*a2*a1*a4)]);;

\bigbreak
\centerline{($\cdots$)}
\bigbreak
\par
\endgroup

\bigbreak
And the central part for $G_2$ is:
\bigbreak
\begingroup
\tt
\parindent=0cm
\obeylines
\bigbreak
\centerline{($\cdots$)}
\bigbreak

camino:=id;
caminom:=id;

camino1:=camino*k3m*k1m;;
camino1m:=k1*k3*caminom;;
sing1:=k2;;
trenza1:=camino1*sing1*camino1m;;
rel:=List(lista,u->igual(u,trenza1));;

camino2:=camino1;;
camino2m:=camino1m;;
sing2:=(k1*k3)\^{ }7*(k2*k1*k3*k2);;
trenza2:=camino2*sing2*camino2m;;
rel:=Union(rel,List(lista,u->igual(u,trenza2)));;

camino3:=camino2*(k1*k3)\^{ }4*(k2*k1*k3*k2);;
camino3m:=(k2m*k3m*k1m*k2m)*(k1m*k3m)\^{ }4*camino2m;;
sing3:=k2;;
trenza3:=camino3*sing3*camino3m;;
rel:=Union(rel,List(lista,u->igual(u,trenza3)));;

camino4:=camino3*k1m*k3m;;
camino4m:=k3*k1*camino3m;;
sing4:=k2\^{ }2;;
trenza4:=camino4*sing4*camino4m;;
rel:=Union(rel,List(lista,u->igual(u,trenza4)));;

rel:=Union(rel,[Image(caminom,a1*a4*a3*a2*a1*a4)]);;

\bigbreak
\centerline{($\cdots$)}
\bigbreak
\par
\endgroup
\finparrafo

\enddocument